\UseRawInputEncoding

\documentclass[12pt]{amsart}
\usepackage{epsfig,color}
\usepackage{CJK}
\usepackage{times}
\usepackage{mathptmx}
\makeatother

\theoremstyle{definition}

\def\fnum{equation} 
\newtheorem{Thm}[\fnum]{Theorem}

\newtheorem{Lem}[\fnum]{Lemma}
\newtheorem{Con}[\fnum]{Conjecture}

\numberwithin{equation}{section}

 \newcommand{\C}{\ensuremath{\mathbb{C}}}

 \newcommand{\ba}{\begin{align*}}
 \newcommand{\ea}{\end{align*}}
 \newcommand{\na}{\nabla}

\newcommand{\p}{\partial}

\title{Review of the Siu-Yang conjecture }

\author{X.Zhao}%
\address{}

\thanks{}

\email{}

\begin{document}

\maketitle

 \begin{abstract}
In this notes we review the development and recent results of the Siu-Yang conjecture which is that every K\"ahler-Einstein compact complex manifold of complex dimension two with negative sectional curvature is biholomorphic to a compact quotient of the complex 2-ball.
  \end{abstract}
  Siu and Yang have proposed a well-known conjecture on compact K\"ahler-Einstein surface on 1981 that \cite{1}: 
\begin{Con}: Every K\"ahler-Einstein compact complex manifold of complex dimension two with negative sectional curvature is biholomorphic to a compact quotient of the complex 2-ball.
\end{Con}
The conjecture is quite important in complex geometry as it characterizes the K\"ahler-Einstein surface with negative sectional curvature, which in some sense is a counterpart of the Frankel conjecture in the negative sectional curvature manifolds. 

It is deserved to be pointed out that the K\"ahler-Einstein condition is not redundant since there exists a counter example constructed by Mostow and Siu that there is indeed a compact K\"ahler surface with negative sectional curvature whose universal covering is not biholomorphic to the complex two-ball, which implies that it is not a quotient of the complex 2-ball. The construction involves finding a subgroup of Aut B generating by three complex reflections and the consideration of singularities when lifting maps and taking quotients. For more details, see \cite{2}

 The conjecture, however, has not been completely proved so far while there are some partial results now. Most results are stated with the assumption that the average of holomorphic sectional curvature is not “too close to” the maximum holomorphic sectional curvature, while the closeness is characterized by the function $\frac{K_{av}-K_{min}}{K_{max}-K_{min}}$. Since the function is trivially less than  $\frac{2}{3}$, whose proof will be repeated in this note, the upper bound would be removable if it is $\frac{2}{3}$. But now the upper bound is $\frac{1}{2}$ which is proved by Daniel Guan \cite{3}, while it is still unknown whether the conjecture is true when the upper bound is larger than $\frac{1}{2}$. Note that Daniel Guan is pessimistic about the conjecture for higher upper bound. 
 
Since Siu and Yang’s paper is the origin of the conjecture and set up the basic notations and ideas for the latter papers, I would like to start from the Siu and Yang’s result and then move on to the Daniel Guan’s result and some other recent developments. 

It is deserved to be pointed out that I have an original observation on Siu and Yang’s proof, which could give a better result than the one stated in that paper.

The reading note is divided into four parts: 

1. Introduction: notations and preliminaries

2. Siu and Yang’s result and further observation

3. D. Guan, D. Chen, Y. Hong and H.C. Yang’s results

4. Other related results and recent developments

\section{introduction}
 Now we denote M as a K\"ahler-Einstein manifold of complex dimension two. Denote P as a point in M. Let $e_1$  and $e_2$ be a unitary frame in the tangent space of M at P and denote $R_{\alpha \bar{\beta}\gamma\bar{\delta}}$ as the components of curvature tensors with respect to the frame $e_1$  and $e_2$, which is the same as the notation on our class. Denote $\rho$ as the scalar curvature at P. Every calculation happens only at a point P, otherwise we will make additional comments on it. $K_{max}(P)$ denote the maximum of holomorphic sectional curvature at P.
 
Since we are in a K\"ahler manifold of complex dimension two, there are only few independent components of curvature tensor, which is definitely good news. In order to further simplify the computation, we would like to fix $e_1$ as a critical direction of the holomorphic sectional curvatures, which will cancel many components of curvature tensor due to the following lemma.

\begin{Lem}If $e_1$ is a critical direction for the holomorphic sectional curvature, then  $R_{\alpha \bar{\beta}\gamma\bar{\delta}}$ vanishes when precisely three of $\alpha,\beta,\gamma,\delta$ are equal.
\label{lem1}\end{Lem}
\begin{proof}  
Let $\zeta=(\zeta^1,\zeta^2)$ on the unit sphere in $\C^2$. Let $\zeta^i=\xi^i+\sqrt{-1}  \eta^i$ for $i=1,2$. Then $\xi^1=\sqrt{1-(\xi^2)^2 -(\eta
^1)^2-(\eta^2)^2}$. Since $e_1$ is a critical direction, we know that the first order derivative of $\sum R_{\alpha\bar{\beta} \gamma\bar{\delta}} \zeta^\alpha \bar{\zeta^\beta} \zeta^\gamma \bar{\zeta^\delta}$ with respect to $\xi^2$ ,$\eta^1$, $\eta^2$ all vanishes at $\zeta=(1,0)$. Then we get 
             \begin{align}
 0=R_{1\bar{1}1\bar{2}}+R_{1\bar{1}2\bar{1}}+R_{1\bar{2}1\bar{1}}+R_{2\bar{1}1\bar{1}}\\
0=-4R_{1\bar{1}1\bar{1}}+4R_{1\bar{1}1\bar{1}}\\
0=\sqrt{-1} (R_{1\bar{1}2\bar{1}}+R_{2\bar{1}1\bar{1}}-R_{1\bar{2}1\bar{1}}-R_{1\bar{1}1\bar{2}})\end{align}
Combine with the property of curvature tensor and K\"ahler-Einstein condition, we get $R_{1\bar{1}1\bar{2}}=\sqrt{-1} R_{1\bar{1}2\bar{1}}=-R_{2\bar{2}1\bar{2}}=-R_{1\bar{1}1\bar{2}}=0$. The lemma is proved. 
\end{proof}

After proving the lemma, we know that there are only 2 independent curvature tensors  $R_{1\bar{1}2\bar{2}}$   and  $R_{1\bar{2}1\bar{2}}$ using the lemma above, K\"ahler-Einstein condition and Bianchi identity. Noted that the Bianchi identity is used to show that
\begin{align*}R_{1\bar{2}2\bar{1}}=-R_{1\bar{1}\bar{2}2}-R_{12\bar{1}\bar{2}}=-R_{1\bar{1}\bar{2}2}=R_{1\bar{1}2\bar{2}}\end{align*}
where $R_{12\bar{1}\bar{2}}=0$ comes from the fact that J is parallel in K\"ahler manifold.

K\"ahler-Einstein condition is used to show that
\begin{align}R_{1\bar{1}1\bar{1}}=\rho-R_{1\bar{1}\bar{2}2}=R_{22\bar{2}\bar{2}}
\end{align}From now we denote,
\begin{align}A=2R_{1\bar{1}2\bar{2}}-R_{1\bar{1}1\bar{1}}\\
B=R_{1\bar{2}1\bar{2}}                   
\end{align}Now we only consider the case that $K_{max} (P)\neq K_{min} (P)$ while the other case will be treated in the second and third part. Let $e_1$ be a critical direction of the holomorphic sectional curvature, in which case we could use the lemma above.

Consider the holomorphic sectional curvature in the direction $\zeta^1 e_1+\zeta^2 e_2$ with $|\zeta^1 |^2+|\zeta^2 |^2=1$, combine with the properties of curvature we proved above, we have

\begin{align}\sum R_{\alpha\bar{\beta}\gamma\bar{\delta}}  \zeta^\alpha \bar{\zeta^\beta}\zeta^\gamma \bar{\zeta^{\delta}}=R_{1\bar{1}1\bar{1}}+2A|\zeta^1 \bar{\zeta^2} |^2+2ReR_{1\bar{2}1\bar{2}} (\zeta^1 \bar{\zeta^2} )^2 .\end{align}
In order to simplify the term $ReR_{1\bar{2}1\bar{2}} (\zeta^1 \bar{\zeta^2} )^2$, let $R_{1\bar{2}1\bar{2}}=|R_{1\bar{2}1\bar{2}}|e^{\sqrt{-1} \theta}$ and $\zeta^1 \bar{\zeta^2} e^{\sqrt{-1}\frac{\theta}{2}}=u+\sqrt{-1}v$, denote
\begin{align}\sigma=A+|B|\\
   \tau=A-|B|  \end{align}
Then we have $\sigma\geq\tau$ and
\begin{align} 
\sum R_{\alpha\bar{\beta}\gamma\bar{\delta}}  \zeta^\alpha \bar{\zeta^\beta} \zeta^\gamma \bar{\zeta^\delta} =R_{1\bar{1}1\bar{1}}+2(\sigma u^2+\tau v^2).\end{align}
From the formula, we know that 
\begin{align}If \,\,R_{1\bar{1}1\bar{1}}=K_{max} (P) \,\,then\,\, \sigma=A+|B|\leq 0\end{align}
\begin{align}If \,\,R_{1\bar{1}1\bar{1}}=K_{min} (P) \,\,then\,\, \tau=A-|B| \geq 0\end{align}
From now on, let $R_{1\bar{1}1\bar{1}}=K_{min} (P)$, then we would like to express the $K_{max} (P)$ and $K_{av}(P)$. Since $\zeta^1 \bar{\zeta^2} e^{\sqrt{-1}  \frac{\theta}{2}}=u+\sqrt{-1} v$, we have
\begin{align}\sum R_{\alpha\bar{\beta}\gamma\bar{\delta}}\zeta^\alpha \bar{\zeta^\beta} \zeta^\gamma \bar{\zeta^\delta}\leq R_{1\bar{1}1\bar{1}}+2\sigma(u^2+v^2 )\leq R_{1\bar{1}1\bar{1}}+\frac{\sigma}{2}\end{align}
Taking $\zeta^1=\frac{1}{\sqrt{2}},\,\zeta^2=\frac{1}{\sqrt{2}}e{\sqrt{-1}\frac{\theta}{2}}$, then we get the equality, thus 
\begin{align}K_{max} (P)=R_{1\bar{1}1\bar{1}}+\frac{\sigma}{2} \end{align}
For $K_{av}(P)$, using the integration 
\begin{align}K_{av} (P)=c\int_{\C^2}(\sum R_{\alpha\bar{\beta}\gamma\bar{\delta}}  \zeta^\alpha \bar{\zeta^\beta} \zeta^\gamma \bar{\zeta^\delta})e^{-r^2}\end{align}
where $c^{-1}=\int_{\C^2}r^4 e^{-r^2}=6 \pi^2$

Since we have the calculus fact that
\begin{align}\int_{\C^2}|\zeta^1|^{2p} |\zeta^2|^{2q} e^{-r^2}=\pi^2 (p!)(q!)\\ 
\int_{\C^2}(\zeta^1)^2 (\bar{\zeta^2} )^2 e^{-r^2}=0 \end{align}
Plug in the $\sum R_{\alpha\bar{\beta}\gamma\bar{\delta}}  \zeta^\alpha \bar{\zeta^\beta} \zeta^\gamma \bar{\zeta^\delta}$, we get 
\begin{align*}\int_{\C^2}(\sum R_{\alpha\bar{\beta}\gamma\bar{\delta}}  \zeta^\alpha \bar{\zeta^\beta} \zeta^\gamma \bar{\zeta^\delta})e^{-r^2} =\pi^2 (2R_{1\bar{1}1\bar{1}}+4R_{1\bar{1}2\bar{2}}+2R_{2\bar{2}2\bar{2}})\end{align*}
Thus, we have
\begin{align}K_{av} (P)=\frac{2}{3}(R_{1\bar{1}1\bar{1}}+R_{1\bar{1}2\bar{2}})\\
K_{av} (P)-K_{min} (P)=\frac{1}{3} A\\
K_{max} (P)-K_{min} (P)=\frac{1}{2}(A+|B|)\end{align}
Since we have made the assumption that $R_{1\bar{1}1\bar{1}} (P)$ attains $K_{min}(P)$, we know that $A\geq|B|$ from previous discussion. Thus 
\begin{align}\frac{1}{3}(K_{max}(P)-K_{min} (P))\leq K_{av} (P)-K_{min} (P)\leq \frac{2}{3}(K_{max}(P)-K_{min} (P))\end{align}
which is mentioned at the beginning.

By now we have only calculated on a single point, which is not good enough since we do not know what happens around the point. It would be reasonable to ask whether there is a local vector field along which the curvature attains maximum or minimum. It turns out that this is indeed the case, which will be proved below.

1) If $\tau=A-|B|>0$, then $K_{min} (P)$ is achieved only at $u=v=0$. i.e. $\zeta^1=0$  or  $\zeta^2=0$. Thus $K_{min}(P)$ is achieved only at two points of $\P(T_{M,P})$ defined by $e_1,e_2$. By the implicit function theorem we can find a smooth unitary frame field $e_1(Q)$,$e_2(Q)$ where Q belongs to a neighborhood of P. 

Moreover, with the assumption $\frac{K_{av}-K_{min}}{K_{max}-K_{min}}<\chi<\frac{2}{3}$, we can get $\psi B\geq A$. Since we have the assumption $K_{max}\neq K_{min}$, thus $B\neq 0$. (Otherwise $A=B=0$,$K_{max}=K_{min}$ ). Since $B=R_{121\bar{2}}\neq 0$, we could define 
\begin{align}\xi(Q)=\frac{1}{\sqrt{2}}(e_1(Q)+(exp(\frac{\sqrt{-1}}{2}arg R_{121\bar{2}(Q)})e_2 (Q))\end{align}

which satisfies that the holomorphic sectional curvature along $\xi(Q)$ is $K_{max}(Q)$

2) If $\tau=A-|B|=0$, then $K_{min} (P)$ is achieved at $u=0$. i.e. $arg \frac{\zeta^2}{\zeta^1} =\pm\frac{\pi}{2}-\frac{\theta}{2}$. This is a circle in $\P(T_{M,P})$. Since $K_{max}\neq K_{min}$, we have $\sigma>0$.
\begin{align} 
\sum R_{\alpha\bar{\beta}\gamma\bar{\delta}}  \zeta^\alpha \bar{\zeta^\beta} \zeta^\gamma \bar{\zeta^\delta}\leq R_{1\bar{1}1\bar{1}}+2\sigma(u^2+v^2 )\leq R_{1\bar{1}1\bar{1}}+\frac{\sigma}{2}\end{align}
The equality holds (i.e. achieve $K_{max}$) if and only if 
\begin{align}0=|\zeta^1|^2+|\zeta^2|^2\pm 2Re(\zeta^1 \bar{\zeta^2} e^{\sqrt{-1}  \frac{\theta}{2}})=|\zeta^1\pm\zeta^2 e^{-\sqrt{-1}  \frac{\theta}{2}}|^2\end{align}
i.e. $\zeta^1=\frac{1}{\sqrt{2}}a$ ,$\zeta^2=\pm \frac{1}{\sqrt{2}} a e^{\sqrt{-1}  \frac{\theta}{2}}$, where a has an absolute value 1. Thus $K_{max}(P)$ is achieved at at only two points in $\P(T_{M,P})$.

To summarize, in the first case there are two smooth vector fields in an open neighborhood of P which attain $K_{max}(Q)$ and $K_{min}(Q)$, while $K_{min}(P)$ is attains at only two points in $\P(T_{M,P})$. Then we achieve the aim we set in the beginning.

In the second case, $K_{min}(P)$ is achieved on a circle and $K_{max}(P)$ is attained at only two points in $\P(T_{M,P})$. It seems that the problem is not easy to be handled in this case. However, if we turns to apply the above argument replacing by $R_{1\bar{1}1\bar{1}}=K_{max}$, we would find that we are in the similar situation of the first since $K_{max}(P)$ is attained at only two points in $\P(T_{M,P})$. From the first case we know that such vector fields also exist.
Overall, there are always two smooth vector fields in an open neighborhood of P which attain $K_{max}(Q)$ and $K_{min}(Q)$.

After getting the formula of $K_{min} (P)$, $K_{av} (P)$ and $K_{max} (P)$, we would like to calculate the covariant derivative and Laplacian of the curvature tensor, which is the main goal for the first part.

From the definition and the fact that g is parallel, we know that 
\begin{align}\sum R_{\alpha\bar{\beta}\gamma\bar{\delta}}  g^{\alpha\bar{\beta}}=\rho g_{\gamma\bar{\delta}}, \,\,\sum\Delta R_{\alpha\bar{\beta}\gamma\bar{\delta}}  g^{\alpha\bar{\beta}}=0\end{align}
where $\rho$ is the scalar curvature.

Since the frame is orthonormal, we only need to calculate $\Delta R_{1\bar{1}1\bar{1}}$ and  $\Delta R_{1\bar{2}1\bar{2}}$. Using Bianchi identity, commutation formula and $K\ddot{a}hler$-Einstein condition, we can get that 
\begin{align*}\delta R_{1\bar{1}1\bar{1}}=\sum_s\na_s \na_{\bar{s}} R_{1\bar{1}1\bar{1}} =\sum_s\na_s \na_{\bar{1}} R_{1\bar{s}1\bar{1}} \\ 
          =\sum_s\na_{\bar{1}} \na_s R_{1\bar{s}1\bar{1}} +\sum_{s,k}R_{s\bar{1}\bar{1}1} R_{k\bar{s}1\bar{1}}+R_{s\bar{1}k\bar{s}} R_{1\bar{k}1\bar{1}} +R_{s\bar{1}\bar{k}1} R_{1\bar{s}k\bar{1}}+R_{s\bar{1}k\bar{1}} R_{1\bar{s}1\bar{k}} \\
          =\sum_s\na_{\bar{1}} \na_1 R_{s\bar{s}1\bar{1}} +\sum_{s,k}R_{s\bar{1}\bar{k}1} R_{k\bar{s}1\bar{1}}+R_{s\bar{1}k\bar{s}} R_{1\bar{k}1\bar{1}} +R_{s\bar{1}\bar{k}1} R_{1\bar{s}k\bar{1}} +R_{s\bar{1}k\bar{1}} R_{1\bar{s}1\bar{k}}\\=-AR_{1\bar{1}2\bar{2}}+|B|^2 \end{align*}
where the last equality is gotten from a straight forward calculation.
Similarly, we can calculate  $\Delta R_{1\bar{2}1\bar{2}}$ and got
  \begin{align}\Delta R_{1\bar{1}1\bar{1}}=-AR_{1\bar{1}2\bar{2}}+|B|^2 \\
\Delta R_{1\bar{2}1\bar{2}}=3(R_{1\bar{1}2\bar{2}}-A)B\end{align}
Recall that we have already proved that there is a local vector field  $\xi(Q)$ along which the holomorphic sectional curvature is equal to $K_{min}(P)$. Choose a local coordinate system $z^1,z^2$ at P such that, 
\begin{align}\xi(P)=\frac{\p}{\p z^1}(P),\eta(P)=\frac{\p}{\p z^2}(P).\end{align}
Denote     $\xi=\sum\xi^i \frac{\p}{\p z^i},\eta=\sum\eta^j \frac{\p}{\p z^j}$
 
Let $S_{\alpha\bar{\beta}\gamma\bar{\delta}}$ be the component of curvature with respect to the basis  $\xi$  and  $\eta$. By replacing $\eta$ by $e^{\sqrt{-1} \theta} \eta$ for some constant  $\theta$, we can assume $S_{1\bar{2}1\bar{2}} (P)$ is real and nonnegative. It is deserved to be pointed out that this type of rotating trick helps us a lot in the below calculation.

Let $g_{\alpha\bar{\beta}} dz^\alpha dz^{\bar{\beta}}$ be the K\"ahler-Einstein metric on M. From $\sum g_{\alpha\bar{\beta}} \xi^\alpha  \bar{\xi^\beta}=1$. Take $\na$, since $\xi(P)=\frac{\p}{\p z^1}(P)$,$\eta(P)=\frac{\p}{\p z^2}(P)$ take value at P we know that $Re \na\xi^1 (P)=0$. Choose a real-valued function $\theta$ such that  $\theta(P)=0$ and  $d\theta(P)=-Im \na\xi^1 (P)$. By replacing $\xi$  by $e^{\sqrt{-1} \theta} \xi$ we can assume that $\na\xi^1 (P)=0$. Similarly, we can assume $\na\eta^2 (P)=0$. Take $\na$ of $\sum g_{\alpha\bar{\beta}} \xi^\alpha  \bar{\eta^\beta}=0$, at P we get 
\begin{align}\na\xi^2 (P)+\na\bar{\eta^1}(P)=0\end{align}
Taking second derivative of $\sum g_{\alpha\bar{\beta}} \xi^\alpha  \bar{\xi^\beta}=1,\sum g_{\alpha\bar{\beta}} \xi^\alpha  \bar{\eta^\beta }=0$  and $\sum g_{\alpha\bar{\beta}} \eta^\alpha  \bar{\eta^\beta}=1$ and taking the value at P, we know that 
\begin{align}
\na_s \na_{\bar{s}} \xi^1 (P)+\na_s \na_{\bar{s}} \bar{\xi^1}(P)+|\na_s \xi^2 |^2 (P)+|\na_{\bar{s}} \xi^2 |^2 (P)=0\\
\na_s \na_{\bar{s}} \xi^2 (P)+\na_s \na_{\bar{s}} \bar{\eta^1}(P)=0\\
\na_s \na_{\bar{s}} \eta^2 (P)+\na_s \na_{\bar{s}} \bar{\eta^2}(P)+|\na_s \eta^1 |^2 (P)+|\na_{\bar{s}} \eta^1 |^2 (P)=0.\end{align}
As for the relationship between $\na R_{\alpha\bar{\beta}\gamma\bar{\delta}}$, $\na\xi$  and $\na\eta$, since
\begin{align}\sum R_{\alpha\bar{\beta}\gamma\bar{\delta}}  \xi^\alpha \bar{\xi^\beta}\xi^\gamma \bar{\eta^\delta}=S_{1\bar{1}1\bar{2}}=0.\end{align}
Take $\na_t$, since we have already chosen the local frame such that the holomorphic sectional curvature attains minimum along the field, we know that $\na S_{1\bar{1}1\bar{2}}=0$. Thus, at P we have 
\begin{align*}\na_t R_{1\bar{1}1\bar{2}}=-A\na_t \xi^2-B\na_t \bar{\xi^2}\end{align*}
Computing the Laplacian of $S_{1\bar{1}1\bar{1}}$, we get that 
\begin{align*}\Delta S_{1\bar{1}1\bar{1}}=\sum\na_s \na_{\bar{s}} \sum R_{\alpha\bar{\beta}\gamma\bar{\delta}}  \xi^\alpha \bar{\xi^\beta}\xi^\gamma  \bar{\xi^\delta}
  =\sum R_{\alpha\bar{\beta}\gamma\bar{\delta}} \na_s \na_{\bar{s}} ( \xi^\alpha \bar{\xi^\beta}\xi^\gamma  \bar{\xi^\delta})+2Re\sum(\na_s R_{\alpha\bar{\beta}\gamma\bar{\delta}})\na_{\bar{s}} ( \xi^\alpha \bar{\xi^\beta}\xi^\gamma \bar{\xi^\delta})+\sum(\na_s \na_{\bar{s}} R_{\alpha\bar{\beta}\gamma\bar{\delta}}) \xi^\alpha \bar{\xi^\beta}\xi^\gamma  \bar{\xi^\delta}.\end{align*}
After straight forward calculation using the bold formula above, and take value at P we get 
\begin{align}\Delta S_{1\bar{1}1\bar{1}} (P)=-2A\sum(|\na_s \xi^2 |^2+|\na_{\bar{s}} \xi^2 |^2)-4B Re\sum(\na_s \xi^2 )(\na_{\bar{s}} \xi^2 )-A R_{1\bar{1}2\bar{2}}+B^2.\end{align}
Similarly, from straight forward calculation we know that 
\begin{align*}Re \Delta S_{1\bar{2}1\bar{2}} (P)=4A Re\sum(\na_s \xi^2)( \na_{\bar{s}} \xi^2)\\+2B\sum(|\na_s \xi^2 |^2+|\na_{\bar{s}} \xi^2 |^2)+3( R_{1\bar{1}2\bar{2}}-A)B\end{align*}
After calculating the Laplacian, we would like to calculate the covariant derivative of S at P. 

\begin{align*}\na_t S_{1\bar{1}1\bar{1}} (P)=\na_t R_{1\bar{1}1\bar{1}} (P)\\
\na_t S_{1\bar{2}1\bar{2}} (P)=\na_t R_{1\bar{2}1\bar{2}} (P)\end{align*}
Since we have proved before that $\na_t R_{1\bar{1}1\bar{2}}=-A\na_t \xi^2-B\na_t \bar{\xi^2}$, we can translate the equation using the second Bianchi identity
\begin{align}
\na_1 R_{1\bar{1}1\bar{1}}=-\na_1 R_{1\bar{1}2\bar{2}}=-\na_2 R_{1\bar{1}1\bar{2}}=A\na_2 \xi^2+B\na_2 \bar{\xi^2}\\
\na_2 R_{1\bar{1}1\bar{1}}=\na_1 R_{1\bar{1}2\bar{1}}=\overline{\na_{\bar{1}} R_{1\bar{1}1\bar{2}}}=-A\na_1 \bar{\xi^2}-B\na_1 \xi^2\\
\na_{\bar{1}} R_{1\bar{2}1\bar{2}}=\na_{\bar{2}} R_{1\bar{1}1\bar{2}}=-A\na_{\bar{2}} \xi^2-B\na_2\bar \bar{\xi^2 }\\
\na_2 R_{1\bar{2}1\bar{2}}=\na_1 R_{1\bar{2}2\bar{2}}=-\na_1 R_{1\bar{1}1\bar{2}}=A\na_1 \xi^2+B\na_1 \bar{\xi^2}.\end{align}

\section{Siu and Yang’s result and some observations}
After calculating the covariant derivative and Laplacian of the curvature tensor, we are ready to give the result of Siu and Yang and our observation.

Siu and Yang’s result is that
\begin{Thm}if M is a compact K\"ahler-Einstein surface with nonpositive holomorphic bisectional curvature and
\begin{align}\frac{K_{av}-K_{min}}{(K_{max}-K_{min}}\leq\chi<\frac{2}{3(1+\sqrt{\frac{6}{11}})}\end{align}
then M is biholomorphic to a compact quotient of the complex 2-ball with an invariant metric or the 2-complex-dimensional plane. 
\end{Thm}
 
Now we give a simplified proof of the theorem and using my own observation to give a better result that the statement is also true, if 
\begin{align*}\frac{K_{av}-K_{min}}{(K_{max}-K_{min}}\leq\chi<\frac{2}{3(1+\sqrt{\frac{1}{6}})}\end{align*} .
At first, we divide M into two parts. N denotes the ball-like points such that $K_{max}=K_{min}$, i.e.
\begin{align*}N=\{p\in M\,|  K_{max} (P)=K_{min} (P)\}\end{align*}
then $M=N\sqcup(M\setminus N)$

We claim that 
       
     If N$\neq$M,then the real dimension of N is $\leq$2

We try to prove the claim by contradiction. If the real dimension if N is 3, then from the $K\ddot{a}hler$-Einstein condition we know that on N
\begin{align}R_{\alpha\bar{\beta}\gamma\bar{\delta}}=\frac{\rho}{3}(g_{\alpha\bar{\beta}} g_{\gamma\bar{\delta}}+g_{\alpha\bar{\delta}}g_{\gamma\bar{\beta}}\end{align}
Thus, the covariant derivative of $R_{\alpha\bar{\beta}\gamma\bar{\delta}}$ along the tangent space of N is 0. If we could show that on N, the covariant derivative of $R_{\alpha\bar{\beta}\gamma\bar{\delta}}$ vanishes along all directions, then combine with the real analyticity of the K\"ahler-Einstein metric we know that M is locally symmetric and $M=N$, which contradicts with the assumption that $M\neq N$. Thus, it suffices to prove that on N, the covariant derivative of $R_{\alpha\bar{\beta}\gamma\bar{\delta}}$ vanishes along all directions.

Denote $e_1$ as the normal vector on N. The proof is a standard induction type argument while we do induction on the number of $e_1$ in $\{e_{v_i} \}$ in $\na_{e_{v_1}}, \na_{e_{v_2}}, \na_{e_{v_3}}...\na_{e_{v_n}}R(e_{\mu_1},e_{\mu_2},e_{\mu_3},e_{\mu_4})$. The key points are:

	1)We have commutation formula
\begin{align*}\na_{e_{v_1}} \na_{e_{v_2}} R(e_{\mu_1 },e_{\mu_2} ,e_{\mu_3},e_{\mu_4})=\na_{e_{v_2}} \na_{e_{v_1}} R(e_{\mu_1},e_{\mu_2 },e_{\mu_3},e_{\mu_4}) .\end{align*}

2) We could strictly reduce the number of  $e_1$ in $\{e_{v_i}\}$ by moving $e_1$ to the last position using 1 and by the Bianchi identity we can move the $e_1$ into R. Here the only concern is that what if the term in R contains  $e_1$ in each pair. This problem could be solved by operate J or using the K\"ahler-Einstein condition. Since the number of  $e_1$ in $\{e_{v_i}\}$ is strictly reduced, we can use the induction hypothesis to get the conclusion. 

After proving that N is M or N’s real dimension is bounded by 2, the second step is to construct a bounded strictly superharmonic function $\Phi$ in $M\setminus N$. 

This implies the theorem since it could be extended to N as the dimension of N is less or equal to 2. Moreover, the extension is still superharmonic and bounded, which must be a constant since the manifold is closed. But it contradicts with the strict superharmonicity of $\Phi$. Thus, the theorem is proved after constructing the function $\Phi$. Indeed, finding such a function is the central work in the proof and also illuminates the following papers on this topic.

Now we construct such function. Define 
\begin{align}\Phi=6|S_{1\bar{2}1\bar{2}} |^2-(S_{1\bar{1}1\bar{1}}-2S_{1\bar{1}2\bar{2}})^2\end{align}
Recall that $S_{1\bar{2}1\bar{2}}(P)$ is real. At P we have $\Phi=6B^2-A^2 $

Then using the formula in 1, we can get 
   \begin{align*}   \Delta\Phi=-12(A^2 -B^2) \sum(|\na_s \xi^2 |^2+|\na_{\bar{s}} \xi^2 |^2 ) + 12B^2 \sum(|\na_s \xi^2 |^2+|\na_{\bar{s}} \xi^2 |^2 )\\   + 24AB Re\sum(\na_s \xi^2 )( \na_{\bar{s}} \xi^2 ) +6(6B^2-A^2 ) R_{1\bar{1}2\bar{2}}\\-30AB^2+6\sum|\na_s R_{1\bar{2}1\bar{2}} |^2 +6\sum|\na_{\bar{s}} R_{1\bar{2}1\bar{2}} |^2\\ 
  -18\sum|\na_s R_{1\bar{1}1\bar{1}} |^2 \end{align*}
Of course, this formula will not satisfy us since it is too complicated. But after using the formulas at the end of part 1 we can transfer $|\na_s R_{\alpha\bar{\beta}\gamma\bar{\delta}} |$ into the form of the first three terms in the formula of $\Delta\Phi$. Also notice that on $M\setminus N$ we have 
\begin{align}C=-6(6B^2-A^2 ) R_{1\bar{1}2\bar{2}}+30AB^2>0.\end{align}
Since M has nonpositive sectional curvature and,
\begin{align}\frac{K_{av}-K_{min}}{(K_{max}-K_{min}}\leq\chi<\frac{2}{3(1+\sqrt{\frac{1}{6}})}\end{align}
which implies that $6B^2-A^2\geq0$.
 
It is deserved to be pointed out that here we only need $\frac{K_{av}-K_{min}}{(K_{max}-K_{min}}\leq\chi<\frac{2}{3(1+\sqrt{\frac{1}{6}})}$ which is weak than the assumption in the theorem. Indeed, the constant in the theorem is needed in a latter step, which actually is not essential using an observation of my own.

     Now we can simplify $\Delta\Phi$ to be
      \begin{align*} \Delta\Phi=-30(A^2 -B^2 )(|\na_2 \xi^2 |^2+|\na_1\bar \xi^2 |^2 )
    -6(A^2 -B^2 )(|\na_1 \xi^2 |^2+|\na_2\bar \xi^2 |^2 )\\-C+ 6|\na_1 R_{1\bar{2}1\bar{2}} |^2\\
                   +6|\na_2\bar R_{1\bar{2}1\bar{2}} |^2.\end{align*}
Note that this is still not always nonpositive. By a regular trick that considering the  $\Phi^\lambda$, we could use $\na_t \Phi$ to make  $\Delta\Phi^\lambda$ nonpositive. 

From the definition of  $\Phi$, we have
\begin{align*}\na_s \Phi=6B(\na_s R_{1\bar{2}1\bar{2}}+\na_s R_{2\bar{1}2\bar{1}})+6AR_{1\bar{1}1\bar{1}}.\end{align*}
By using the formula of  $\na_s R$ we can simplify the formula to be

\begin{align*}\na_1 \Phi=6B\na_1 R_{1\bar{2}1\bar{2}}+6(A^2-B^2)\na_2 \xi^2,\end{align*}

\begin{align*}\na_2 \Phi=6B\na_2 R_{2\bar{1}2\bar{1}}-6(A^2-B^2)\na_1 \bar{\xi^2}.\end{align*}
To make  $\Delta\Phi^\lambda>0$, More precisely to make 
 \begin{align*} \Delta\Phi^\lambda=\lambda\Phi^(\lambda-2) (\Phi\Delta\Phi-(1-\lambda) \sum|\na_s \Phi|^2 )>0.\end{align*}
Since we have $\Phi>0$ and $C>0$, where C appears in $\Delta\Phi$ and guarantee the strict inequality. Thus, to get the inequality we only need
\begin{align*}(6B^2-A^2 )6|\na_1 R_{1\bar{2}1\bar{2}} |^2\leq(1-\lambda)36|B\na_1 R_{1\bar{2}1\bar{2}}+(A^2-B^2 ) \na_2 \xi^2 |^2\\
                                 +(6B^2-A^2 )30(A^2-B^2)|\na_2 \xi^2 |^2,\end{align*}
And, 
\begin{align*}(6B^2-A^2 )6|\na_2\bar R_{1\bar{2}1\bar{2}} |^2\leq(1-\lambda)36|B\na_2\bar R_{1\bar{2}1\bar{2}}-(A^2-B^2 ) \na_1\bar \xi^2 |^2\\  +(6B^2-A^2 )30(A^2-B^2)|\na_1\bar \xi^2 |^2.\end{align*}
Since the two inequalities are similar, we only need to prove the first one. The inequality is equivalent to
\begin{align*}|B\na_1 R_{1\bar{2}1\bar{2}} |^2\leq\frac{6(1-\lambda) B^2}{6B^2-A^2} |B\na_1 R_{1\bar{2}1\bar{2}}+(A^2-B^2 ) \na_2 \xi^2 |^2\\
+\frac{5B^2}{A^2-B^2 } |(A^2-B^2)\na_2 \xi^2 |^2\end{align*}
Denote $f=B\na_1 R_{1\bar{2}1\bar{2}}+(A^2-B^2 ) \na_2 \xi^2,\,g=-(A^2-B^2)\na_2 \xi^2$, we only need to show
\begin{align*}|f+g|^2\leq\frac{6(1-\lambda) B^2}{6B^2-A^2}|f|^2+\frac{5B^2}{A^2-B^2}|g|^2\end{align*}
i.e.
\begin{align*}2|fg|\leq(\frac{6(1-\lambda) B^2}{6B^2-A^2}-1)|f|^2+(\frac{5B^2}{A^2-B^2}-1)|g|^2.\end{align*}
From the AM-GM inequality, we only need to show that 
\begin{align*}1\leq(\frac{6(1-\lambda) B^2}{6B^2-A^2}-1)(\frac{5B^2}{A^2-B^2}-1).\end{align*}
At this point, in Siu and Yang’s paper they use the assumption
$\frac{K_{av}-K_{min}}{(K_{max}-K_{min}}<\frac{2}{3(1+\sqrt{\frac{6}{11}})}$ to deal with the inequality which is strong than $\frac{K_{av}-K_{min}}{(K_{max}-K_{min}}<\frac{2}{3(1+\sqrt{\frac{1}{6}})}$. It is deserved to be pointed out that the stronger assumption is only used here. Thus, if we could find a better way to deal with the inequality, we can get a better result than Siu and Yang. Now indeed this could be achieved by an original observation that the left-hand side is equal to 
\begin{align*}\frac{A^2-6\lambda B^2}{A^2-B^2}\end{align*}
Since we have made the assumption that $R_{1\bar{1}1\bar{1}}(P)$ is the minimal sectional curvature at P, we know that $A\geq B\geq 0$. With the weaker condition  $\frac{K_{av}-K_{min}}{(K_{max}-K_{min}}<\frac{2}{3(1+\sqrt{\frac{1}{6}})}$, we know that $6B^2-A^2\geq 0$, thus the AM-GM could be used since they are all nonnegative. If we take  $\lambda\leq \frac{1}{6}$, then
\begin{align*}(\frac{6(1-\lambda) B^2}{6B^2-A^2}-1)(\frac{5B^2}{A^2-B^2 }-1)=\frac{A^2-6\lambda B^2}{A^2-B^2}\leq 1.\end{align*}
Thus the inequality is proved. So we have
  \begin{align*}\Delta\Phi^\lambda>0 \,\,on\,\, M\setminus N\end{align*}
As we said previously, this implies the theorem since it could be extended to N as the dimension of N is less or equal to 2. Moreover, the extension is still superharmonic and bounded, which must be a constant since the manifold is closed. But it contradicts with the strict superharmonicity of $\Phi^\lambda$. Thus $N=M$ and M has constant holomorphic sectional curvature, which leads to the conclusion of the theorem.

By now the result in Yang and Siu’s paper is completely stated and actually we improve the result from $\frac{2}{3(1+\sqrt{\frac{6}{11}})}$ to a better constant $\frac{2}{3(1+\sqrt{\frac{1}{6}})}$ using the observation above.

3. D. Guan, D. Chen, Y. Hong and H.C. Yang’s results
After giving the result of Siu and Yang, we are now ready to state the most recent result of this problem published on 2017 \cite{3}.
\begin{Thm}
if M is a compact K\"ahler-Einstein surface with nonpositive holomorphic bisectional curvature and
\begin{align*}\frac{K_{av}-K_{min}}{(K_{max}-K_{min}}\leq \frac{1}{2}\end{align*}
Then M is biholomorphic to either the compact quotient of the complex 2-ball with an invariant metric or the 2-complex-dimensional plane\end{Thm}

It is deserved to point out that here there is another possibility that M is the 2-complex-dimensional plane. For more general case the conjecture may not be true.

The proof is divided into several steps. Here we give the outline of the proof and details can be found in Guan’s paper \cite{3}.

Firstly, we show that under the condition of the theorem, there must be a point at which $A=B=0$. i.e. $N\neq \emptyset$
 
The key point of the first step is to consider the function
\begin{align*}\Phi_1=\frac{|B|^2}{A^2}\end{align*}
By calculating the $\Delta\Phi_1$ and consider the minimum point $x_0$ of $\Phi_1$, at which $\Delta\Phi_1\geq 0$, we know that if $x_0\notin N$, i.e. $A\neq 0$ or $B\neq 0$ , then $A\equiv B$ on M. This is a straight forward calculation combining with simple linear algebra. Once we have $A\equiv B$ on M, the assumption in Siu and Yang’s theorem is satisfied. Thus we can use Siu and Yang’s result to conclude that $N=M$ which implies that $N\neq \emptyset$.

In particular we know that if  $\frac{K_{av}-K_{min}}{K_{max}-K_{min}}(P)=\frac{1}{2}$, then $P\in N$. Thus if  $3B(P)=A(P)$, we know that $P\in N$. So $3B-A\neq 0$  on  N

Secondly, we need to consider a different function which is given by Hong Cang Yang \cite{4}
\begin{align*}f=(3B-A)^{\frac{1}{3}} \,\,   on\,\,   M\setminus N\end{align*}
Here $\frac{1}{3}$ is carefully chosen to make sure that $\Delta f<0$ on $M\setminus N$. and $3B-A$ is from the assumption that $\frac{K_{av}-K_{min}}{(K_{max}-K_{min}}\leq \frac{1}{2}$, which is equal to $3B-A\geq 0$. The calculation is straight forward once we fix the strategy. 

Since $\Delta f<0$ on $M\setminus N$, we can use the same argument in Siu and Yang’s result to conclude that $N=M$. So M has constant nonnegative holomorphic sectional curvature, which implies that M is a complex quotient of the complex 2-ball with an invariant metric or the 2-complex-dimensional plane.

\section{Other related results and recent developments}
After giving the result of Guan, which is a landmark result of this problem, I would also like to state some other results in different settings or with more assumptions. 

Fangyang Zheng \cite{5} has proved that for a compact $K\ddot{a}hler$-Einstein surface with nonpositive bisectional curvature and $c_1<0$. If we further assume that M is not quasi-ample, then it is a quotient of bidisc $D\times D$. Where $\Omega$ is quasi-ample if the tautological line bundle L is ample and $Y\cdot L^{dimY}>0$ for any irreducible subvariety $Y\subset P(T_M)$ with $\pi(Y)=M$ where $\pi:P(T_M)\to M$ is the projection.

Vestislav Apostolov and Johann Davidov \cite{6} have classified compact Hermitian surface of nonnegative isotropic curvature. More precisely we can either say that such manifold is biholomorphically isometric to one of four types of manifolds or know that it is biholomorphic to one of two types of manifolds. 

Atreyee Bhattacharya and Harish Seshadri \cite{7} have proved that for a compact Ricci-flat 4-manifold, if $K_{max} (p)\leq-cK_{min} (p)$ for all $p\in M$. Where $0\leq c<\frac{2+\sqrt6}{4}$, then $(M,g)$ is flat. Also similarly for Ricci-flat K\"ahler surfaces if $K_{max} (p)\leq-cK_{min} (p)$ for all $p\in M$, where $0\leq c<\frac{1+\sqrt{3}}{2}$ then $(M,g)$ is flat.

At the end of the reading note, we would like to give a perspective of the conjecture. The basic idea is to prove that the surface has constant holomorphic sectional curvature. As pointed out in the Guan’s paper, improving by modifying the upper bound of   $\frac{K_{av}-K_{min}}{(K_{max}-K_{min}}$ is quite hard now since the test function is difficult to find using the trick of taking power when we increase the upper bound to more than $\frac{1}{2}$. If we still want to walk along this path, finding new ways of finding test function is crucial.

\end{document}